# A Theory for Coloring Walks in a Digraph

*Seth Chaiken*

Department of Computer Science, LI 67A

State University of New York at Albany

Albany, New York 12222

sdc@cs.albany.edu

*ABSTRACT*

Consider edge colorings of directed graphs where edges of the form $v_1 v_2$ and $v_2 v_3$ must have different colors. Here, $v_1 \neq v_2$, $v_2 \neq v_3$ but $v_1 = v_3$ is possible. It is known that this coloring induces a vertex coloring by sets of edge colors, in which edge $v_1 v_2$ in the graph implies that the set color of $v_1$ contains an element not in the set color of $v_2$; conversely, each such set coloring of vertices induces one or more edge colorings. We show that these relationships generalize to colorings of of k(vertex)-walks in which two k-walks have different colors if one is the prefix and the other is the suffix of a common (k+1)-walk. For full generality the colors belong to a partially ordered set P; and the prefix color $c_1$ and the suffix color $c_2$ must satisfy $c_1 \nleq c_2$. The set color construction generalizes to generating the lower order ideal in P from a set of k-walk colors; these order ideals (antichains in P, equivalently) are partially ordered by containment. We conclude that a P coloring of k-walks exists if and only if there is a vertex coloring by $A^{k-1}(P)$, where A is the operator that maps a poset to its poset of lower order ideals, due to Birkhoff. In the case when the graph G is symmetric, this condition means that the largest antichain size (Dilworth number) of $A^{k-1}(P)$ is at least the chromatic number $\chi(G)$. Thus the directed chromatic index problem is generalized and reduced to poset coloring of vertices.

This work uses ideas, results and motivations due to Cole and Vishkin [Deterministic coin tossing with applications to optimal parallel list ranking, Information and Control 70 (1986), pp. 32-53] and to Becker and Simon [How robust is the n-cube?, Information and Commputation 77 (1988), pp. 162-178] on vertex covers for subsets of (n-2)-cubes. Some relationships to Ramsey theory are sketched.

## 1. Introduction

In this paper all graphs, partially ordered sets, and other sets will be finite unless otherwise stated.

Let $G = (V,E)$ be a digraph with vertices V and edges E. Each edge is an ordered pair (i,j) of distinct vertices; that is, we assume no multiple edges or self loops in G. We denote (i,j) is an edge in G by $i \rightarrow j$. For $k \geq 1$, a k-walk is a sequence $(v_1 v_2 ... v_k)$ of vertices for which $v_1 \rightarrow v_2 \rightarrow ... \rightarrow v_k$. The length of a k-walk is k-1. The edges of the walk are $\{(v_i, v_{i+1}) \mid 1 \leq i < k\}$.

Let S be a set. A coloring of the k-walks of G is a function c from k-walks into S such that for every (k+1)-walk $(v_1 v_2 ... v_k v_{k+1})$,

$$c(v_1 v_2 ... v_k) \neq c(v_2 v_3 ... v_{k+1}).$$

We will see that for $k > 1$, a necessary and sufficient condition for the existence of this coloring is that there exists a function $c'$ on the (k-1)-walks into $2^S$ (the collection of subsets of S) such that for every k-walk $(v_1 v_2 ... v_k)$,

$$c'(v_1 v_2 ... v_{k-1}) \nsubseteq c'(v_2 v_3 ... v_k).$$

For this reason we define at the onset the notion of P-coloring of k-walks, where P is a partially ordered set (poset).

**Definition:** A P-coloring of k-walks of G is a mapping $c_k$ from k-walks in G to P such that for every (k+1)-walk $(v_1 v_2 ... v_{k+1})$ in G,

$$c_k(v_1 v_2 ... v_k) \nleq c_k(v_2 v_3 ... v_{k+1}).$$

This paper generalizes ideas and results due to Cole and Vishkin [7] on their deterministic symmetry breaking computation step used in parallel algorithms; and Becker and Simon [2] on how many node faults in Boolean n-cube networks certain fault tolerant parallel processing schemes can tolerate. Becker and Simon defined the "directed chromatic index" of a directed graph and studied it using vertex colors in one to one correspondance with sets of edge colors.

In both subjects, an edge coloring $c_2$ by elements of S is related a vertex coloring $c_1$ by subsets of S so that the color $x = c_2(v_1, v_2)$ of edge $v_1 v_2$ satisfies $x \in c_1(v_1)$ and $x \notin c_1(v_2)$. Given an edge coloring, a vertex coloring is given by $c_1(v_1) = \{ c_2(v_1 v_2) \mid v_1 \to v_2 \}$. Conversely, given a vertex coloring, an edge coloring $c_2$ is obtained by setting $c_2(v_1 v_2)$ to any element of $c_1(v_1) - c_1(v_2)$, whenever $v_1 \to v_2$. In our generalization, S is a poset P. The formation $\{ c_2(v_1 v_2) \mid v_1 \to v_2 \}$ becomes the lower order ideal generated by these $\{c_2(v_1 v_2)\}$. The fundamental construction of the distributive lattice of order ideals A(P) (equivalently, antichains) of a finite partial order P was given by Birkhoff [1,3]. It is used in the theorem that every finite distributive lattice is represented by A(P) where P is the poset of join irreducibles $\neq 0$ in the lattice. Therefore, the vertex color domain $2^S$ generalizes to this lattice of order ideals or antichains. We also find that the ''reduction'' process from an edge P-coloring to a vertex A(P)-coloring naturally extends to reducing a coloring of k-walks to a coloring of (k-1)-walks, for any $k > 1$. Since A(P) is a poset, reduction can be iterated until a vertex coloring is obtained. In applications, lower bounds on P for which a P-coloring of k-walks exists might be derived from lower bounds for ordinary vertex coloring of the undirected version of the given graph.

The motivating problem from Becker and Simon [7] is to find the minimum number of processor nodes in a Boolean n-cube network whose failure guarantees that a given testing scheme will not find at least one n-2 dimensional subcube all of whose processors are good. In the model for this problem, each node has coordinates $b_1 b_2 \cdots b_n$, $b_i \in \{0,1\}$, $V = \{1,2,...,n\}$ and $i \to j$ means that the (n-2)-subcube with $b_i = 0$ and $b_j = 1$ will be tested. Suppose the edges are S-colored. All the tests corresponding to the edges with any one color can be made to fail with one faulty processor. Thus the desired minimum number of processors is the minimum $|S|$ for which an S-coloring exists. Becker and Simon call this the directed chromatic index of the graph. They use Sperner's theorem to show that the directed chromatic index of the complete symmetric directed graph on n vertices is the smallest r for which $\binom{r}{\lfloor r/2 \rfloor} \geq n$. This r is the smallest $|S|$ for which $2^S$ contains an antichain of n sets. Our result that relates P-coloring of edges to A(P)-coloring of vertices generalizes their work that relates directed chromatic index to chromatic number.

''Deterministic coin-tossing'' or symmetry breaking [7,9] is a way of implementing the recursive pairing technique [13] efficiently for parallel algorithms [12,14] when random operations are not permitted. We have a list of n distinct nodes, each having its own processor. Each processor uses the same program on different data (the SIMD model). In one time step, a processor may only access data in its own node and in its neighboring nodes. We would like each processor to simultaneously calculate a color for its node from a small fixed set of colors so that adjacent nodes receive different colors. That way, an independent set of nodes that are spaced no more than a constant distance apart within the list can be identified; this is called a ruling set in [7]. A ruling set determines a set of non-adjacent edges which when contracted, reduces the list length by at least a constant fraction. Thus the list can be contracted to a single node in O(log n) parallel contraction steps. It is assumed that each node initially has a data item, such as an address or processor number, that is distinct from that of any other node. There is no assumed relationship between these initial data and the order of the nodes in the list. (However, the processors for the first and last nodes do ''know'' they are first or last respectively.)

The deterministic symmetry breaking step of [7] can be seen in our terms. For simplicity, assume that a processor may only access data in the node (if any) to its right and in its own node. We begin with an n vertex simple directed path in which each vertex initially has a unique color. For each path edge $v_i v_{i+1}$ in parallel, the algorithm repeats the following steps:

(1) Determine $c_2(v_i v_{i+1})$ from the colors $c_1(v_i)$ and $c_1(v_{i+1})$.

(2) Assign the new color $c_1(v_i) := c_2(v_i v_{i+1})$.

The domain of the color sets changes with each repetition. The domain size decreases logarithmically with each repetition until it becomes a constant which is determined by algorithm details. Let $\log^* n$ denotes the number of times to iterate the $\log_2$ function starting on n to get a value $\leq 2$. In $O(\log^* n)$ repetitions, the algorithm calculates a vertex coloring with a constant number of colors, assuming that the initial color domain has size O(n). Further refinements [9] can produce a 3-coloring--these refinements require processing colors of 3 adjacent vertices for the last steps.

The entire algorithm determines a color for each vertex $v_i$ (except for the last $O(\log^* n)$ of them) from the initial color sequence of vertices $v_i, v_{i+1}, \cdots v_{i+l}$ where $l=O(\log^* n)$. Contiguous subsequences of algorithm steps calculate newer colors from older colors similarly, except that the old color sequences may have non-adjacent repeated colors. Thus, every contiguous subsequence of steps calculates a coloring of some k-walks in a directed graph; the vertex set is the domain of the older colors. The graph generally is complete and symmetric because the algorithm must work on all initial color sequences.

Examples and applications of P-coloring of k-walks in graph theory are given in section 2. Sections 3 and 4 contain the main results. Section 5 explains some simple consequences and how the work of Becker and Simon [7] on directed chromatic index is generalized. The paper concludes with section 6 which mentions connections with Ramsey theory and other literature on edge colorings, and two open questions.

## 2. Examples

In this section we give examples of P-colorings and other problems in graph theory that can be formulated in terms of them.

(1) A coloring of k-walks into S is a special case P-coloring because we can take P=S with the trivial partial order $x \leq y$ if and only if $x = y$.

(2) Ordinary vertex coloring of an undirected graph is a special case of (1), where $k = 1$ and G contains $v_1 v_2$ and $v_2 v_1$ for each undirected edge $\{v_1, v_2\}$.

(3) The directed chromatic index of Becker and Simon[2] in a directed graph G is the minimum size of a finite set S such that the edges, i.e. the 2-walks, of G are P-colorable with trivial $P = S$. Example (6) below is a generalization of this.

(4) The undirected version of a directed graph G has the same vertex set and contains edge $\{v_1, v_2\}$ if and only if $v_1 \rightarrow v_2$ or $v_2 \rightarrow v_1$. $\chi(G)$ denotes the chromatic number of the undirected version of G. A P-coloring of the 1-walks of G is *a fortiori* an ordinary vertex coloring of the undirected version of G, hence $|P| \geq \chi(G)$ if the 1-walks (i.e., vertices) of G can be P-colored.

Becker and Simon [2] introduce the notion of a "simplest directed version" $G_0$ of an undirected graph G, which is an orientation of G that minimizes the directed chromatic index. They prove that this minimum $c'(G_0) = \lceil \log_2 \chi(G) \rceil$. The following remark will help us generalize their result: If $|P| \geq \chi(G)$ then G can be oriented so its vertices can be P-colored. For, any classical coloring c of the vertices by elements of P induces orientations so $u \rightarrow v$ implies $c(u) \not\leq c(v)$. Furthermore, it is easy to use a topological order of P to induce a P colorable acyclic orientation.

(5) If $P = C_n$, the chain of n elements, then the vertices of G are P-colorable if and only if there is no k-walk in G with $k > n$. In particular, it is necessary for G to be acyclic.

(6) Let $S = \{1,2,...,n\}$ be a set of colors and for each color i, let $l_i$ be an integer satisfying $l_i \geq 1$. We are interested in coloring the edges (i.e., 2-walks) of G so there is no walk of length $l_i + 1$ all of whose edges have color i. (This condition implies that there is no monochromatic directed cycle.) Let P be the product partial order on $[0,l_1] \times [0,l_2] \times ... \times [0,l_n]$; that is, $(x_1, x_2, ... x_n) \leq (y_1, y_2, ... y_n)$ if and only if $x_i \leq y_i$ for all i, $1 \leq i \leq n$. The i-th component of $x \in P$ is denoted $x(i)$.

**Theorem.** G can be edge colored as above if and only if the 1-walks (i.e., the vertices) of G can be P-colored.

**Proof.** (only if) Let d be a coloring of the edges that satisfies the above conditions. Let $c: V \rightarrow P$ be defined by $c(v)(i)$=maximum length of a walk starting from vertex v all of whose edges are colored i. Suppose $v_1 \rightarrow v_2$ and $d(v_1, v_2) = i$. Then $c(v_1)(i) > c(v_2)(i)$; hence $c(v_1) \not\leq c(v_2)$.

(if) For any edge $v_1 \rightarrow v_2$, $c(v_1) \not\leq c(v_2)$ so there exists an i with $c(v_1)(i) > c(v_2)(i)$. Let $d(v_1, v_2)$ be any such i. Since by hypothesis $0 \leq c(v_1)(i) \leq l_i$, there cannot be a walk starting at $v_1$ of length exceeding $l_i$ all of whose edges have d color i.

A necessary condition that G can be edge colored as above is that $\chi(G) \leq |P| = \Pi(l_i+1)$. A result of Busolini [5] can now be derived from $|V|/\alpha \leq \chi(G)$ where $\alpha$ is the independence number of G. We should note that the necessity of $\chi(G) \leq \Pi(l_i+1)$ is implied by Chvátal's result [6] which states that this is necessary for the edges to be colored so there is no *simple* path of length $(l_i+1)$ all of whose edges have color i.

We can further generalize this example for $k \geq 2$ so that we color the k-walks rather than single edges. The condition that there is no walk whose length exceeds $l_i$ and all of whose edges are colored i becomes:

(*) For all walks $(v_1 v_2 ... v_k ... v_{k+l-1})$ in G,

(**) $c(v_1 v_2 ... v_k) = c(v_2 v_3 ... v_{k+1}) = \cdots = c(v_l v_{l+1} ... v_{k+l-1}) = i$

implies $l \leq l_i$.

The existence of such a coloring is equivalent to there being a P-coloring of the (k-1)-walks of G. Again, $P = [0, l_1] \times [0, l_2] \times ... \times [0, l_n]$.

To prove this, suppose first that $c_k$ is an $S = \{1, 2, ..., n\}$ coloring of k-walks that satisfies (*). We claim $c_{k-1}: (k-1)\text{-walks} \to P$ is a P-coloring, where

$c_{k-1}(v_1 v_2 ... v_{k-1})(i) = \max \{l \mid \text{there exists a walk } (v_k v_{k+1} ... v_{k+l-1}), \text{ such that (**)}\}$.

Now $0 \leq c_{k-1}(v_1 ... v_{k-1})(i) \leq l_i$. Suppose $(v_1 v_2 ... v_k)$ is a k-walk. Let $i = c_k(v_1 v_2 ... v_k)$. Then $c_{k-1}(v_1 v_2 ... v_{k-1})(i) > c_{k-1}(v_2 v_3 ... v_k)(i)$ because for any $l$ and $v_{k+1}, v_{k+2} ... v_{k+l}$ for which

$c_k(v_2 ... v_{k+1}) = c(v_3 ... v_{k+2}) = \cdots = c_k(v_{l+1} ... v_{l+k}) = i$,

also $c_k(v_1 v_2 ... v_k) = i$, so $c_{k-1}(v_1 v_2 ... v_{k-1})(i) \geq l+1$.

Conversely, suppose the (k-1)-walks are P-colored. For each k-walk choose $c_k(v_1 v_2 ... v_k)$ to be any i for which $c_{k-1}(v_1 v_2 ... v_{k-1})(i) > c_{k-1}(v_2 v_3 ... v_k)(i)$. Then if there is any $l$ and walk $(v_1 v_2 ... v_k ... v_{k+l-1})$ for which (**) is true, then

$l_i \geq c_{k-1}(v_1 ... v_{k-1})(i) > c_{k-1}(v_2 ... v_k)(i) > \cdots > c_{k-1}(v_l ... v_{k+l-1})(i) > c_{k-1}(v_{l+1} ... v_{k+l-1})(i) \geq 0$

so $l \leq l_i$.

## 3. The Reduction Theorem

The main result is that the k-walks of G are P-colorable if and only if the (k-1) walks of G are $A(P)$-colorable, where $A(P)$ is a poset on the antichains of P (actually $A(P)$ is a distributive lattice) given by Birkhoff [1,3], see also Greene and Kleitman [10]. Recall that an antichain of P is a subset of P in which every pair of distinct elements are incomparable. We will use the fact that $A(P)$ is order isomorphic to the collection of the order ideals of P.

**Definition.** Let $A(P)$ denote the collection of all antichains on P. For $X, Y \in A(P)$, let $X \leq Y$ if and only if for every element $x \in X$ there is a $y \in Y$ so $x \leq y$.

**Remark.** When P is the trivial partial order on a set S, $A(P)$ coincides with $2^S$, with $X \leq Y$ if and only if $X \subseteq Y$.

**Definition.** If $X \subseteq P$ let MAX(X) denote the subset of maximal elements in X; that is $x \in $ MAX(X) if and only if $x \in X$ and there is no $y \in X$ with $x < y$. Note that MAX(S) is an antichain.

**Definition** If $X \subseteq P$, X is called a (lower) order ideal if whenever $y \in X$ and $x \leq y$ then $x \in X$. I(Y) denotes the smallest order ideal of those that contain Y. Note that $x \in I(Y)$ if and only if there exists $y \in Y$ where $x \leq y$.

The function $X \to $ MAX(X) is a one-to-one correspondence from the order ideals of P onto $A(P)$. Its inverse is $Y \to I(Y)$. In fact, $I_1$ and $I_2$ are order ideals with $I_1 \subseteq I_2$ if and only if MAX($I_1$) $\leq$ MAX($I_2$) (as antichains). Since the collection of order ideals is closed under union and intersection, it comprises a distributive lattice of subsets under $\subseteq$. Hence $A(P)$ is a distributive lattice. Birkhoff's theorem states that every finite distributive lattice L is isomorphic to $A(P)$ where $P = J(L) \subseteq L$ is the poset of the join-irreducibles $\neq 0$ in L. See [1,3].

**Theorem.** Let $k > 1$. The k-walks of G are P-colorable if and only if the (k-1)-walks of G are $A(P)$ colorable.

**Proof.** (only if) Suppose $c_k$ is a P-coloring of the k-walks of G. We claim that an $A(P)$-coloring of the k-1 walks of G, denoted by $c_{k-1}$, is

$$c_{k-1}(v_1v_2...v_{k-1}) = \text{MAX}(\{c_k(v_1v_2...v_{k-1}v_k) \mid (v_1v_2...v_k) \text{ is a } k\text{-walk}\}).$$

To verify the claim, suppose $(v_1v_2...v_k)$ is a k-walk. Let

$$X = I(c_{k-1}(v_1v_2...v_{k-1})) = I(\{c_k(v_1v_2...v_{k-1}v'_k) \mid (v_1v_2...v_{k-1}v'_k) \text{ is a } k\text{-walk}\}).$$
$$Y = I(c_{k-1}(v_2v_3...v_k)) = I(\{c_k(v_2v_3...v_kv'_{k+1}) \mid (v_2v_3...v_kv'_{k+1}) \text{ is a } k\text{-walk}\}).$$

We assert that $x = c_k(v_1v_2...v_k) \in X$ and $x \notin Y$, so $X \not\subseteq Y$ and therefore $c_{k-1}(v_1v_2...v_{k-1}) \not\leq c_{k-1}(v_2v_3...v_k)$. First, $x \in X$ since $(v_1v_2...v_k)$ is a k-walk. Suppose now, contrary to the claim, that $x \in Y$. There must then be a k-walk $(v_2v_3...v_kv'_{k+1})$ and P-color $z = c_k(v_2v_3...v_kv'_{k+1})$ such that $x \leq z$. However the relationship $x \leq z$ contradicts the assumption that $c_k$ is a P-coloring because $(v_1v_2...v_kv'_{k+1})$ is a (k+1)-walk.

(if) Suppose $c_{k-1}$ is an $A$(P)-coloring of the (k-1)-walks of G. For each k-walk $(v_1v_2...v_k)$, there exists $x \in P$ such that

$$x \in c_{k-1}(v_1v_2...v_{k-1}) - I(c_{k-1}(v_2v_3...v_k)).$$

To prove this, we note that the condition that $c_{k-1}$ is an $A$(P)-coloring means that

$$c_{k-1}(v_1v_2...v_{k-1}) \not\leq c_{k-1}(v_2v_3...v_k)$$

so

$$I(c_{k-1}(v_1v_2...v_{k-1})) - I(c_{k-1}(v_2v_3...v_k)) \neq \emptyset$$

contains some $z \in P$ and so an $x \in c_{k-1}(v_1v_2...v_{k-1})$ with $x \geq z$ exists. $x \notin I(c_{k-1}(v_2v_3...v_k))$ because if not, $z \in I(c_{k-1}(v_2v_3...v_k))$ as well. Let us take for $c_k(v_1v_2...v_k)$ any such $x \in P$. To show that $c_k$ is a P-coloring of k-walks, suppose $(v_1v_2...v_kv_{k+1})$ is a (k+1)-walk and the P-color chosen for $c_k(v_2v_3...v_kv_{k+1})$ is

$$y \in c_{k-1}(v_2v_3...v_k) - I(c_{k-1}(v_3v_4...v_{k+1})).$$

Suppose $x \leq y$. $I(c_{k-1}(v_2v_3...v_k))$ is an order ideal that contains y so it must then contain x. This contradicts the choice of x. Hence $x \not\leq y$ and $c_k$ is a P-coloring of k-walks.

**Remark.** The set from which $c_k(v_1v_2...v_k)$ is chosen may be enlarged to

$$I(c_{k-1}(v_1v_2...v_{k-1})) - I(c_{k-1}(v_2v_3...v_k)).$$

## 4. Expansion Theorems

It is obvious that if $c_k$ is a P-coloring of the k-walks of G, then $c_{k+1}$ defined by $c_{k+1}(v_1 v_2 ... v_k v_{k+1}) = c_k(v_1 v_2 ... v_k)$ is a P-coloring of the (k+1) walks of G. Thus

**Trivial Expansion Theorem:** If the k-walks of G are P-colorable, then so are the k′−walks of G, for any $k' \geq k$.

In all the examples in Section 1 with non-trivial P, P was in fact a distributive lattice. The theorem of Birkhoff characterizes distributive lattices as those lattices that are isomorphic to a lattice of order ideals of a poset under set union and intersection. The representing poset is order isomorphic to the set of join irreducibles J(P) other than 0 in the lattice P. This is one example of the representation r of a poset P by a collection of subsets so $x \leq y$ in P if and only if $r(x) \subseteq r(y)$. In general, whenever P has such a representation, a P-coloring of k-walks can be expanded to (k+1)-walks as follows.

**Expansion Theorem:** Suppose P has a representation $r: P \to 2^S$ so $x \leq y$ in P if and only if $r(x) \subseteq r(y)$. If the k-walks of G are P-colorable, then the (k+1)-walks of G are S-colorable where the partial order on S is trivial.

**Proof.** Suppose c is the P-coloring. Let $(v_1 v_2 ... v_{k+1})$ be any k+1-walk. Hence $c(v_1 v_2 ... v_k) \not\leq c(v_2 v_3 ... v_{k+1})$ so $r(c(v_1 v_2 ... v_k)) \not\subseteq r(c(v_2 v_3 ... v_{k+1}))$. Let the color $c'(v_1 v_2 ... v_{k+1})$ be any $i \in r(c(v_1 ... v_k)) \setminus r(c(v_2 ... v_{k+1}))$. We conclude $i \neq c'(v_2 v_3 ... v_{k+2})$ because $c'(v_2 v_3 ... v_{k+2}) \in r(c(v_2 v_3 ... v_{k+1}))$.

**Example.** Suppose the vertices of G are n-colored by the trivial poset $P = \{c_1, c_2, ... c_n\}$. Suppose $n \leq \begin{bmatrix} r \\ \lceil r/2 \rceil \end{bmatrix}$. Then P is representable by any n distinct subsets of size $\lceil \frac{r}{2} \rceil$ of $S = \{1, 2, ..., r\}$. It follows that the edges of G are S-colorable.

When P is a distributive lattice, a P-coloring of the k-walks of G can be expanded to J(P)-coloring of the (k+1)-walks. That is, in condition $c'(v_1 v_2 ... v_{k+1}) \neq c'(v_2 v_3 ... v_{k+2})$, the $\neq$ is strengthened to $\not\leq$. For in the above proof, if $i \leq c'(v_2 v_3 ... v_{k+2})$ then $i \in r(c(v_2 v_3 ... v_{k+1}))$ because $r(c(v_2 v_3 ... v_{k+1}))$ is an order ideal. (Compare to the ''if'' part of the reduction theorem.)

## 5. Conditions for P-colorings

Recall that the minimum element 0 of poset P, if it exists, is the element such that $0 \leq x$ for all $x \in P$. A minimal element a has the property that there is no $x \in P$ for which $x < a$. Analogous definitions hold for maximum and maximal elements. Every finite poset has at least one minimal element. A minimal element a is unique if and only if P has a minimum element, which would then be a.

Suppose c is a P-coloring of the k-walks of G. If P has a (unique) minimum element 0, then $c(v_1 v_2 ... v_k) = 0$ implies $v_k$ is a sink or is isolated. Conversely, if the k-walks of G can be P-colored, then there exists a P-coloring c for which $c(v_1 v_2 ... v_k)$ is a minimal element whenever $v_k$ is a sink. An analogous result holds for a maximum element, sources, and maximal elements with $v_k$ replaced by $v_1$.

**Example.** Let P be the poset $\{ a, b, 0, 1 \}$ with $0 < a < 1$, $0 < b < 1$ and a,b incomparable. The vertices of G are P-colorable if and only if a P-coloring can be found in which every source is colored 1, every sink is colored 0 (isolated vertices can be of any color), and every other vertex is colored a or b. Thus G is P-colorable if and only if G′ is bipartite where G′ is obtained from G by deleting every source and sink.

The complete symmetric digraph $K_n$ on n vertices has $u \to v$ and $v \to u$ for every pair of distinct vertices u and v.

Generally speaking, the larger k is the smaller $|P|$ is needed for the k-walks of a given graph to be P-colored. Indeed, the constructions in [9] demonstrate that the k-walks of the complete symmetric directed graph with n vertices are 3-colorable provided that $k = \Omega(\log^* n)$. If G is acyclic and k−1 exceeds the length L of the longest walk in G, the k-walks are $\phi$-colorable. In that case, the reduction theorem's construction colors each (L+1)-walk with $\phi$. However, the k-walks of some cyclic graphs cannot be P-colored if $|P| \leq 2$ no matter how large k is. The proof is left to the reader.

**Theorem.** Suppose G has an odd directed cycle $v_1, v_2, ..., v_m$. If $|P| < 3$, no P-coloring of the k-walks exists, for any $k \geq 1$.

Becker and Simon [2] provided results on three problems concerning the coloring of edges. Their results relate the minimum number of colors required to the ordinary chromatic number $\chi(G)$. The problems, stated in our terminology, together with our generalizations, are as follows:

(1) ''Directed chromatic index problem''.

Let G be a directed graph. The chromatic index of G, c(G) is the least integer n such that the 2-walks, i.e., edges of G can be S-colored where S is the trivially ordered poset with n elements. Becker and Simon show $c(G) \geq \lceil \log_2 \chi(G) \rceil$. (Note $\lceil \log_2 m \rceil = \min\{k \mid m \leq 2^k\}$.) This is equivalent to $\chi(G) \leq 2^n$ is necessary for the edges of G to be n-colored. The immediate generalization is that $\chi(G) \leq |A(P)|$ is necessary for the edges to be P-colored. This results from the reduction theorem and the remark that $\chi(G) \leq |A(P)|$ is necessary for the vertices to be $A(P)$-colored.

(2) ''Chromatic index of the doubly directed version of $G_0$.''

Let $G_0$ be the undirected version of directed graph G. $c''(G_0)$ is the chromatic index of $G_2$, where digraph $G_2$ has the same vertices $G_0$ has and for every edge (u,v) in $G_0$, G has $u \to v$ and $v \to u$. Note that $c(G) \leq c''(G_0)$. Becker and Simon show $c''(G_0) = R_\chi(G)$ where $R_m = \min\{k \mid m \leq \binom{k}{\lceil k/2 \rceil}\}$ Their proof is immediately generalized to demonstrate that the vertices of $G_0$ are $A(P)$–colorable so adjacent vertices have incomparable colors if and only if $\chi(G_0)$ is no more than the maximum size of an antichain of $A(P)$. (N.B. This is an ''antichain of antichains''.) Let the maximum size of a antichain of poset $A$ be denoted by Dil($A$). Thus the edges of $G_2$ are P-colorable if and only if $\chi(G) \leq \text{Dil}(A(P))$.

(3) ''Chromatic index of some simplest directed version of $G_0$.''

Let $G_0$ be an undirected graph. $c'(G_0)$ is the chromatic index of an oriented version $G'$ of $G_0$ with least chromatic index. An oriented version of $G_0$ is obtained by replacing each edge (u,v) in $G_0$ by exactly one of the directed edges $u \to v$ and $v \to u$. Becker and Simon show $c'(G_0) = \lceil \log_2 \chi(G) \rceil$. Note $c'(G_0) \leq c(G)$. Our generalization is that

$\chi(G_0) \leq |A(P)|$ is necessary and sufficient for there to exist an orientation $G'$ of $G_0$ that can be P-colored. Furthermore such an orientation can always be found that is acyclic. The proof combines the reduction theorem with the remark under example (4) of section 2, applied to an $A(P)$-coloring of the vertices. Thus, any G without 2-cycles can be reoriented so that the bound given by (1) becomes exact.

With the next proposition, we can characterize the directed chromatic index of complete tournaments. Let length(G) denote the length (number of edges) of a longest walk in G; length(G) = ∞ if G has a directed circuit. The following statement is vacuous unless G is acyclic.

**Proposition.** If length(G) + 1 ≤ $|A(P)|$ then the 2-walks of G can be P-colored.

**Proof.** Let $A_0 A_1 A_2 ...$ be a reverse topological ordering of $A(P)$ so that $i < j$ implies $A_i \nsubseteq A_j$. If the length(G) + 1 ≤ ($A(P)$) then the vertices of G put into at most $|A(P)|$ classes so class(v)=i when the longest path ending with v has length i. No two adjacent vertices are in the same class. Furthermore, if $u \rightarrow v$ then class(u) < class(v). Therefore the assignment of $A_i$ to each vertex in class i is an $A(P)$–coloring of the vertices. The reduction theorem now implies that the 2-walks of G are P-colorable.

**Example.** Take P to be trivially ordered with $|P| = k$. Then $|A(P)| = 2^k$. So if $\lceil \log_2(\text{length}(G) + 1) \rceil \leq k$ then $c(G) \leq k$. Together with (1) we conclude

$$\lceil \log_2 \chi(G) \rceil \leq c(G) \leq \lceil \log_2(\text{length}(G) + 1) \rceil.$$

When G is the complete tournament on n vertices, $\chi(G) = n = \text{length}(G) + 1$ so the directed chromatic index is $\lceil \log_2 n \rceil$.

## 6. Comments

We thank Professor Stefan Burr answering some of our queries, giving some relationships to other points of view and providing some references. First, problems such as finding the minimum number of colors sufficient to color all the k-walks of the complete digraph on n vertices can be formulated as Ramsey-type questions. Second, the problem of whether the 2-walks of a given digraph can be colored with q colors is a special case of the problem: Given positive integers $l_1, l_2, ..., l_q$ is there a coloring of the 2-walks so there is no walk of length $l_i$ edges all of whose edges have color i? Burr explained how a technique of [5] shows that $n \leq l_1 l_2 ... l_q$ is necessary (and sufficient) for the edges of the transitive tournament on n vertices to be so colored [11]. In general, $\chi(G) \leq l_1 l_2 ... l_q$ is necessary where $\chi(G)$ is the vertex chromatic number of the arbitrary digraph G [5],[6]. ([11] and [6] actually use [8] to prove stronger results where the conditions on walks are limited to paths.) Suppose the edges of G are so colored. A vertex coloring by q-tuples $x_v = (x_1, x_2, ..., x_q)$ with $0 \leq x_i < l_i$ is then induced by $x_v =$ (length of a longest walk of edges colored i starting at vertex v). (See also [4] for application of this idea.) For the transitive tournament each pair of vertices u,v must have $x_u \neq x_v$. In general, if $u \rightarrow v$ then $x_u \not\leq x_v$ in the product partial order.

We close with two open questions: (1) Suppose, in the definition of P-coloring, the sequences $(v_1 v_2 \cdots v_k)$, $(v_2 v_3 \cdots v_{k+1})$ and $(v_1 v_2 \cdots v_k v_{k+1})$ range over a given set of sequences instead of walks in a digraph. Characterize those sets of sequences for which the analog of the reduction theorem of section 2 is true. (2) Can results such as those in [6] which concern coloring of paths rather than walks and which depend on Gallai's theorem [8] that $\chi(G) <$ (the length of any directed path in G) be related to the theory of P-colorings?